\documentclass{article}
\usepackage[T2A]{fontenc}
\usepackage[english]{babel}
\usepackage[tbtags]{amsmath}
\usepackage{amsfonts,amssymb,mathrsfs,amscd,amsmath,comment}

\begin{document}

\author{B.\,S.~Safin\footnote{Nizhnij Novgorod, Russia; E-mail: {\tt boris-nn@yandex.ru}}}
\title {New relations for the old triangle}

\maketitle
\begin{abstract}
In this note we show that in addition to two integers forming a Py\-tha\-go\-rean triple, there also exist two
irrational numbers in terms of which this Pythagorean triple can also be obtained.
We also put forward a~relation between these two pairs and the hypotenuse of a Pythagorean triangle.
\end{abstract}

Assume that we are given a Pythagorean triangle   $a^2+b^2=c^2$ with mutually prime $a, b$, where $a$ is even, and $b$~is odd.
This triangle  is known to be generated by 2~numbers, say integers $(x_1 , y_1$). How to find $ x_1 ,y_1 $ from the available numbers $a,b,c$?

Numbers $x_1 ,y_1$ can be found from the formulas:
\begin{equation}\label{eq1}
x_1=\sqrt {c-b/2}, \quad y_1=\sqrt {c+b/2}.
\end{equation}
These formulas are obtained by simple substitution into formulas for Pythagorean triples:
$$
 a=2 x_1 y_1, \quad b= x_1^2-y_1^2, \quad c= x_1^2+y_1^2.
 $$
It turns out, that in addition to two rational numbers there are two irrational numbers that form the same Pythagorean triple  $a,b,c$.

These numbers $(x_2 ,y_2)$ are calculated as follows:
\begin{equation}\label{eq2}
x_2=\sqrt {(c-a)/2},\qquad  y_2=\sqrt {(c+a)/2}.
\end{equation}

\smallskip

{\bf Example:} $3^2+4^2=5^2$. Here, $x_1=\sqrt {(5-3)/2}=1$, $y_1=\sqrt {(5+3)/2}=2$,  $ x_2=\sqrt{( 5-4)/2}=1/\sqrt 2$,  $y_2=\sqrt {(5+4)/2}=3/\sqrt 2$.
Integer numbers are obtained when taking difference or sum between the hypotenuse and the odd leg. Irrational numbers are obtained
when taking difference between the hypotenuse and the even leg.

Whereas everything is unambiguous for a primitive Pythagorean triangle (that is, we known \textit{a~priori} what number (integer or irrational)
we need to find), two variants are possible for a~triangle in which legs and the hypotenuse have a~common divisor. Here either
both pairs $(x_3,y_3)$, $( x_4,y_4)$ may be irrational (example: $9^2+12^2=15^2$) or one pair is integer and the second one is irrational.
 However, to determine which is integer and which is irrational we can employ only the method of exhaustion (example: $6^2+8^2=10^2$). For
 a~Pythagorean triple with common divisor, there should be two new pairs (say, $(x_3,y_3)$, $(x_4,y_4)$), which generate this triangle.

All this is easily visualized by setting
\begin{equation}\label{eq3}
a=\beta m,\ \ b=\beta n,\ \ c=\beta k \qquad \text{or} \qquad  (\beta m)^2+(\beta n)^2=(\beta k)^2.                                                                                               \end{equation}

Now was assume that the equation $a^2+b^2=c^2$ has common divisor~$\beta$.
Substituting  \eqref{eq3} into \eqref{eq1} and \eqref{eq2} and applying a~little algebra,
we obtain new expressions for 2~pairs of numbers that form a~Pythagorean triple:
\begin{gather}
(\beta m)^2+(\beta n)^2=(\beta k)^2=a^2+b^2=c^2, \nonumber \\
x_3=\sqrt {\beta /2(x_1-y_1)}, \qquad y_3=\sqrt {\beta /2(x_1+,y_1)}, \label{eq4} \\
x_4= x_1\sqrt {\beta},\qquad  y_4= y_1\sqrt \beta .                                                                                                                                        \label{eq5}
\end{gather}

Now let us examine what relations exist between integer and irrational pairs. The analysis that follows applies to any triangles,
as easily follows from direct substitution of formulas \eqref{eq4},\eqref{eq5} into formulas for Pythagorean triples. Hence, we shall
 assume that $ (x_1, y_1)$ are integer numbers, and ($x_2 ,y_2)$ are irrational numbers.

As a result,
\begin{gather*}
2x_1y_1= x_2^2-y_2^2, \\
2x_2y_2= x_1^2-y_1^2, \\
x_1^2+y_1^2= x_2^2+y_2^2,
\end{gather*}
and so the Pythagorean theorem can be rewritten as follows:
$$
 4x_1^2y_1^2+4x_2^2y_2^2= (x_1^2-y_1^2)^2+(x_2^2-y_2^2)^2 = (x_1^2+y_1^2)^2=(x_2^2+y_2^2)^2.
$$
All these relations are proved by direct substitution of formulas \eqref{eq1} and~\eqref{eq2}.

Now let us consider the expression:
$$
 (x_1+y_1)^4=x^4_1+4x^3_1y_1+6x^2_1y^2_1+4x_1y^3_1+y^4_1= x^4_1+ 2x^2_1y^2_1+ y^4_1+ 4x^2_1y^2_1+
4x_1y_1(x_1^2+y_1^2).
$$
Since  $x_1^2+y_1^2= x_2^2+y_2^2$, $ 4x_1y_1= 2(x_2^2-y_2^2)$, $4x^2_1y^2_1=( x_2^2-y_2^2)^2$,
the previous expression equals
$$
 (x_2^2+y_2^2)^2+(x_2^2-y_2^2)^2+2(x_2^4-y_2^4).
$$
Expanding and rearranging the resulting expression, we have
\begin{equation}
 (x_1+y_1)^4=4x_2^4.                                                                                                                                                  \label{eq6}
\end{equation}
A similar argument gives
\begin{equation}
(x_1-y_1)^4=4y_2^4.
 \label{eq7}
\end{equation}

Let us find the sum of the fourth powers of both pairs for the generators of the Pythagorean triangle.
To do so we take formulas \eqref{eq1}, \eqref{eq2}, raise each to the fourth power and add them up. The resulting
expression is
\begin {equation}
2x_1^4+2y_1^4+2x_2^4+2y_2^4=3g^2,
 \label{eq8}
\end{equation}
where $g$ is the hypotenuse of any Pythagorean triangle .

\smallskip

To summarize:

1. Any Pythagorean triangle  has 2~pairs of generators, which are found by formulas \eqref{eq1} and~\eqref{eq2}.

2. For a primitive Pythagorean triangle, one pair is necessary composed of integer numbers, and the other one, of irrational numbers.

3 In Pythagorean triple is nonprimitive, these pairs may be both irrational, or one pair is integer, and the other one, irrational.

4. The pairs of generators $(x_1,y_1)$, $( x_2, y_2)$ for Pythagorean triangle are related as follows:
\begin{align*}
 &  4x_2^4= (x_1+y_1)^4, \\
 &  4y_2^4= (x_1-y_1)^4, \\
 &  2x_1y_1= x^2_2-y^2_2, \\
 & 2x_2y_2= x_1^2-y_1^2, \\
 & x_1^2+y_1^2= x_2^2+y_2^2  .
   \end{align*}

5. The hypotenuse and pairs of numbers are related by the formula $3g^2=2x_1^4+2y_1^4+2x_2^4+2y_2^4$.

6. The Pythagorean theorem can be rewritten as follows:
$$
 4x_1^2y_1^2+4x_2^2y_2^2= (x_1^2-y_1^2)^2+(x_2^2-y_2^2)^2 = (x_1^2+y_1^2)^2=(x_2^2+y_2^2)^2.
$$

\end{document}